\documentclass[12pt]{article}

\input epsf
\usepackage{amssymb}
\newcounter{bla}

\def\Fgauss#1#2#3#4{
{}_2F_1\left(
\begin{array}{c}
\begin{array}{cc} \hskip-10pt#1,{\ } #2 \end{array}\\
\begin{array}{c} \hskip-10pt#3 \end{array}
\end{array}
\hskip-8pt;\,#4
\right)}

\def\tfrac#1#2{{{\lower.6ex
\hbox{$\scriptstyle#1$}}\over
{\raise.7ex
\hbox{$\scriptstyle#2$}}}}

\def\arccosh{{\rm arccosh}}
\def\arccot{{\rm arccot}}
\def\dsp#1{\displaystyle#1}

\def\Frac#1#2{\frac{\displaystyle{#1}}{\displaystyle{#2}}}

\begin{document}

\title{An improved algorithm and a Fortran 90 module for computing the conical function $P^m_{-1/2+i\tau}(x)$}

\author{
    Amparo Gil
       \\
       Depto. de Matem\'atica Aplicada y Ciencias de la Comput. \\
       Universidad de Cantabria. 39005-Santander, Spain\\
     { \small e-mail: {\tt amparo.gil@unican.es}}
\\
\\
    Javier Segura
       \\
        Depto. de Matem\'aticas, Estad\'{\i}stica y Comput. \\
        Universidad de Cantabria. 39005-Santander, Spain\\
     { \small e-mail: {\tt javier.segura@unican.es}}
\\
\\
    Nico M. Temme\footnote{Emeritus researcher at Centrum Wiskunde \& Informatica (CWI), 
        Science Park 123, 1098 XG Amsterdam,  The Netherlands}
       \\
        IAA, Abcoude 1391 VD 18,   The Netherlands\\
     { \small e-mail: {\tt Nico.Temme@cwi.nl}}
}

\date{\ }

\maketitle
\begin{abstract}

In this paper we describe an algorithm and a Fortran 90 module ({\bf Conical}) for the
computation of the conical function 
$P^m_{-\tfrac12+i\tau}(x)$ for $x>-1$, $m \ge 0$, $\tau >0$.
These functions appear in 
the solution of Dirichlet problems for domains
bounded by cones; because of this, they are involved in a large number of applications in Engineering and Physics. 

In the Fortran 90 module, the admissible parameter ranges for computing the conical functions in standard IEEE double
precision arithmetic 
are restricted  to $(x,m,\tau) \in (-1,1) \times
[0,\,40] \times [0,\,100]$ and $(x,m,\tau) \in (1,100) \times
[0,\,100] \times [0,\,100]$. Based on tests of the three-term recurrence relation satisfied by these functions
and direct comparison with Maple, we claim a  relative accuracy close to $10^{-12}$
in the full parameter range, although a mild loss of accuracy can be found at some points of the oscillatory region of the
conical functions. The relative accuracy increases to $10^{-13}\,-\,10^{-14}$ in the region of the monotonic
regime of the functions where integral representations
are computed ($-1<x<0$).

\end{abstract}

\section{Introduction}

Conical functions \cite{Duns:2010:Con} (also called Mehler functions) appear in a large 
number of applications in engineering, applied physics
\cite{thebault:2004:geo},
\cite{passian:2010:nano}, particle physics (related to the amplitude for Yukawa potential 
scattering) or cosmology \cite{aplic2:1995:NVL}, among others.    
However, as far as the authors know, the only existing code for computing conical 
functions is given by K\"olbig \cite{Kolbig:1981:PCF}, which is restricted 
for $m=0,1$ (i.e. the functions $P^0_{-\tfrac12+i\tau}(x)$ and
 $P^1_{-\tfrac12+i\tau}(x)$).

In this paper we describe an algorithm and a Fortran 90 module for the
computation of the conical function 
$P^m_{-\tfrac12+i\tau}(x)$ for $x>-1$, $m \ge 0$, $\tau >0$.
The algorithm is based on the use of different methods of computation, depending on the range
of the parameters:
quadrature methods, recurrence relations and  uniform 
asymptotic expansions in terms of elementary functions or in terms of 
modified Bessel function $K_{ia}(x)$ and its derivative $K'_{ia}(x)$. 

The suggested algorithm in \cite{gil:2009:con}
is improved by considering an additional asymptotic expansion
for large $\tau$, which enables to enlarge the range of computation in the $\tau$ variable.   
Also, the algorithm makes use of an expansion in terms of elementary functions in the oscillatory regime,
which was not previously considered in \cite{gil:2009:con}.

Based on direct comparison with Maple and tests of three-term recurrence relations
satisfied by the functions, we claim a relative
accuracy close to $10^{-12}$ (for IEEE standard double precision arithmetic) in the admissible range of parameters for conical functions
in the module {\bf Conical}: $(x,m,\tau) \in (-1,1) \times
[0,\,40] \times [0,\,100]$ and $(x,m,\tau) \in (1,100) \times
[0,\,100] \times [0,\,100]$.

\section{Theoretical background}

  Conical functions $P^m_{-\frac12+i\tau}(x)$ are solutions of the associated Legendre equation

\begin{equation} 
(1-x^2)\Frac{d^2 w}{dx^2}-2x\Frac{dw}{dx}+\left(\nu(\nu+1)-\Frac{m^2}{1-x^2}\right)w=0
\end{equation}
for $\nu=-\frac12+i\tau$ and $x>-1$, $\tau >0$ and $m=0,1,2,\ldots$

The conical function $P^m_{-\frac12+i\tau}(x)$  can be written in terms of the Gauss 
hypergeometric function $_2F_1$ as:

\begin{equation}
\label{changem}
\begin{array}{lcl}
P^m_{-\frac12+i\tau}(x)&=&\cosh(\pi\tau)\Frac{|\Gamma (m+1/2+i\tau)|^2}{\pi \Gamma(1+m)}
\left|\Frac{1-x}{1+x}\right|^{m/2} \times \\ 
&&\Fgauss{\tfrac12 -i \tau}{\tfrac12 +i \tau }{1+m}{\tfrac12-\tfrac12x}.
\end{array}
\end{equation}

The absolute value $\left|\Frac{1-x}{1+x}\right|^{m/2}$ in the
previous formula is the standard normalization which gives real values
for all $x>-1$. 

\begin{figure}
\caption{Graph of the function $P^5_{-\frac12+i}(x)$.
\label{Fig1}}
\begin{center}
\epsfxsize=13cm \epsfbox{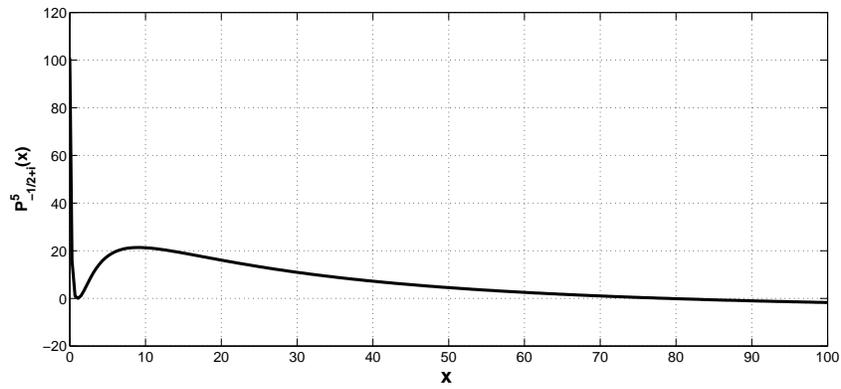}
\end{center}
\end{figure}

\begin{figure}
\caption{Graph of the function $P^5_{-\frac12+i100}(x)$.
\label{Fig2}}
\begin{center}
\epsfxsize=13cm \epsfbox{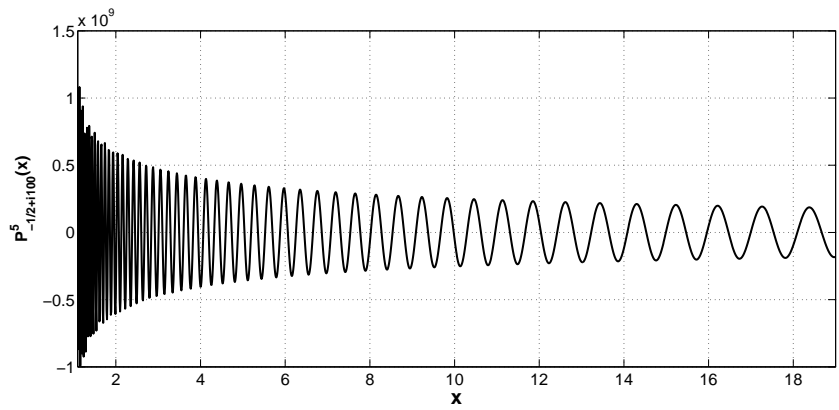}
\end{center}
\end{figure}

The conical functions are monotonic in the interval $(-1,x_c)$ and oscillating in $(x_c,+\infty)$, where
$x_c=\sqrt{1+\beta^2}/\beta$ and $\beta=\tau/m$. In the oscillatory region, the functions strongly oscillate as $\tau$
is taken large. This is apparent in Figures \ref{Fig1} and \ref{Fig2}, where a plot of the functions $P^5_{-\frac12+i}(x)$ and  
$P^{5}_{-\frac12+i100}(x)$, respectively, is shown. 
Figure \ref{Fig2} also shows that the frequency of oscillations is higher for small $x$. 

 Next, we are going to describe the theoretical expressions involved in the 
computation of conical functions:

\subsection{Computation of $P^m_{-\tfrac12+i\tau}(x)$ for $x>0$} 

Two kind of asymptotic expansions are considered: for large $m$ and for large~$\tau$. 

\subsubsection{Asymptotic expansions for large $m$}

For large values of the parameter $m$, asymptotic expansions for $0<x<1$ and $x > 1$, respectively, are used:

\begin{enumerate}
\item The following asymptotic expansion is valid for $0<x<1$, large positive values of $m$ and uniformly
valid for $\tau \ge 0$:

\begin{equation}\label{Pas}
P_{-\frac12+i\tau}^{m}(x)\sim
 \sqrt{\frac{p}{x m}}\,
 \frac{\Gamma(\frac12+m)\,(1-x^2)^{m/2}\cosh(\pi \tau)e^{-m\phi(t_0)}}{\pi}
\sum_{k=0}^\infty \frac{u_{k}(\beta,p)}{m^k}.
\end{equation}

The quantities $\beta$,  $p$ and $\phi(t_0)$ are given by
\begin{equation}\label{betap2}
\beta=\frac{\tau}{m},\quad p=\frac{x}{\sqrt{1+\beta^2(1-x^2)}},
\end{equation}
and 
\begin{equation}\label{phit02}
\phi(t_0)=\ln\frac{x(p+1)}{p(\beta^2+1)}+\beta \arccos \frac{x(1-p\beta^2)}{p(1+\beta^2)}.
\end{equation}
The first few coefficients of the expansion in (\ref{Pas}) are
\begin{equation}\label{uk01}
\begin{array}{l}
u_0(\beta,p)=1,\quad  u_1(\beta,p)=-\dsp{\frac{-\beta^2+5\beta^2p^3-3\beta^2p+3p}{24(\beta^2+1)}},
\\ \\
u_2(\beta,p)=\dsp{\frac{1}{1152(\beta^2+1)^2}}[385\beta^4p^6+462\beta^2(1-\beta^2)p^4-10\beta^4p^3\\ \\
\quad\quad
+(81\beta^4-522\beta^2+81)p^2+6\beta^2(\beta^2-1)p+\beta^4+72\beta^2-72].
\end{array}
\end{equation}
 
\item For $x>1$, we  use a representation in terms of the modified Bessel function $K_{i\tau}(m \zeta)$, which is valid 
for $m$ positive:

\begin{equation}\label{PzwK}
\begin{array}{l}
P_{-\frac12+i\tau}^{m}(x) =
\dsp{\frac{2\Gamma(\tfrac12+m)\,(x^2-1)^{m/2}\,\cosh(\pi \tau)e^{-m \lambda}}{\pi\sqrt{2\pi}\,}}\ \Phi(\zeta) \times \\  \\
\quad\quad\quad
\left[A_m(\beta,\zeta)K_{i\tau}(m\zeta)-B_m(\beta,\zeta)K_{i\tau}^\prime(m\zeta)\right],
\end{array}
\end{equation}
where 
\begin{equation}\label{lambda}
\lambda=\frac12\left(\ln\frac{x^2-1}{\beta^2+1}+\beta\arccos\frac{1-\beta^2}{1+\beta^2}\right),
\end{equation}
\begin{equation}\label{Phizeta}
\beta=\frac{\tau}{m},\quad\Phi(\zeta)=\left(\frac{\zeta^2-\beta^2}{1+\beta^2(1-x^2)}\right)^{\frac14},
\end{equation}
and
the functions $A_m(\beta,\zeta)$ and $B_m(\beta,\zeta)$ have the expansions
\begin{equation}\label{AmuBmu}
A_m(\beta,\zeta)\sim\sum_{n=0}^\infty\frac{A_n(\beta,\zeta)}{m^n},\quad
B_m(\beta,\zeta)\sim\sum_{n=0}^\infty\frac{B_n(\beta,\zeta)}{m^n}.
\end{equation}

These expressions are valid for large $m$.

In this representation, the parameter $\zeta$ and the coefficients of the expansions are given
in two different $x$-regions: $1< x \le x_c$ (the monotonic interval) and $x> x_c$ (the oscillatory region), where
$x_c=\Frac{\sqrt{1+\beta^2}}{\beta}$:

\begin{description}
\item{Case $1< x \le x_c$:}

 In this case the quantity $\zeta\ge \beta$ is given by the implicit equation
\begin{equation}\label{lhrh}
2\left[\sqrt{\zeta^2-\beta^2}-\beta\arccos(\beta/\zeta)\right]=\ln\frac{p+1}{p-1}-\beta\arccos\frac{\beta^2p^2-1}{\beta^2p^2+1},
\end{equation}
where $p$ is given by
\begin{equation}\label{pz}
p=\frac{x}{\sqrt{1+\beta^2(1-x^2)}}.
\end{equation}
This implicit equation cannot be inverted analitically. Then, a method
for computing numerical approximations to the solution of this equation is needed.
In the algorithm, we choose Newton's method given that initial approximations which
guarantee convergence of the method can be obtained \cite{gil:2009:con}.  

The first few coefficients $A_n(\beta,\zeta), B_n(\beta,\zeta)$ in (\ref{AmuBmu}) are
\begin{equation}\label{AB001}
A_0(\beta,\zeta)=1,  \quad B_0(\beta,\zeta)=0, \quad
A_1(\beta,\zeta)=\frac{\beta^2}{24(1+\beta^2)},  
\end{equation}
\begin{equation}\label{B1}
B_1(\beta,\zeta)=-\frac{(5\beta^2(W^3p^3-1-\beta^2)+3W^2(Wp(1-\beta^2)-1-\beta^2)\zeta}{24W^4(1+\beta^2)},
\end{equation}
where $p$ is given in (\ref{pz}) and 
\begin{equation}\label{Wzeta}
W= \sqrt{\zeta^2-\beta^2}.
\end{equation}

\item{Case $x \ge x_c$:}

The quantity $\zeta\in[0, \beta]$ is given by the implicit equation
\begin{equation}\label{lhrh2}
2\left[\sqrt{\beta^2-\zeta^2}-\beta \arccosh(\beta/\zeta)\right]=2\,\arccot\,q-\beta\ln\frac{\beta q+1}{\beta q-1},
\end{equation}
where $q$ is given by
\begin{equation}\label{qz}
q=\frac{x}{\sqrt{\beta^2(x^2-1)-1}}.
\end{equation}

As before, this equation is solved by using Newton's method with appropriated
starting values \cite{gil:2009:con}.

The coefficients $A_0(\beta,\zeta), B_0(\beta,\zeta), A_1(\beta,\zeta)$ are as in (\ref{AB001}), whereas the coefficient $B_1(\beta,\zeta)$ is given by
\begin{equation}\label{B12}
B_1(\beta,\zeta)=-\frac{(5\beta^2(V^3q^3-1-\beta^2)+3V^2(Vq(1-\beta^2)-1-\beta^2)\zeta}{24V^4(1+\beta^2)},
\end{equation}
where $V= \sqrt{\beta^2-\zeta^2}$. 

\end{description}

For the oscillatory case and far away from the transition point
between monotonic and oscillatory behaviour of the function ($x>>x_c>1$), we also use an asymptotic expansion
(valid for large $m$) in 
terms of elementary functions:

\begin{equation}\label{exposc12}
\begin{array}{ll}
\dsp{P_{-\frac12+i\tau}^{m}(x)\sim
2\sqrt{\frac{q}{m x}}\frac{\left(\beta^2+1\right)^{\mu/2}\Gamma\left(\frac12+m\right)}
{\pi}
\cosh (\pi \tau)e^{-\tau(\pi-\arccot(\beta))}}\ \times\\[8pt]
\quad\quad\quad\quad\quad\quad\quad\quad
\dsp{\left(
\cos\chi\sum_{k=0}^\infty\frac{v_k}{m^k}-
\sin\chi\sum_{k=0}^\infty\frac{w_k}{m^k}
\right)}.
\end{array}
\end{equation}
where 

\begin{equation}\label{exposc10a}
\chi=\mu(\beta\xi-\arccot\,q)-\tfrac14\pi,
\end{equation}
and the other parameters used in the expansion are given by

\begin{equation}\label{exposc10b}
\beta=\frac{\tau}{\mu},\quad x_c=\frac{\sqrt{1+\beta^2}}{\beta},\quad
\xi=\arccosh\frac{x}{x_c},\quad q=\frac{\cosh \xi}{\beta\sinh \xi}.
\end{equation}

The first coefficients of the expansion are
{\small
\begin{equation}\label{exposc13}
\begin{array}{ll}
v_0= 1,\quad w_0= 0,\\[8pt]
\dsp{v_1= \frac{\beta^2}{24(1+\beta^2)},}  \\[8pt]
\dsp{w_1=- \frac{q(5\beta^2q^2+3\beta^2-3)}{24(1+\beta^2)}, }  \\[8pt]
\dsp{v_2= -\frac{385\beta^4q^6+462\beta^2(\beta^2-1)q^4+(81\beta^4-522\beta^2+81)q^2-\beta^4-72\beta^2
		+72}{1152(1+\beta^2)^2},}   \\[8pt]
\dsp{w_2= \frac{\beta^2q(5\beta^2q^2+3\beta^2-3)}{576(1+\beta^2)^2}.}
\end{array}
\end{equation}
}
These coefficients also follow from the coefficents $u_k$ given in (4.4) of \cite{gil:2009:con}
 by writing $p=iq$. Then 
\begin{equation}\label{exposc14}
v_k=(-1)^k\Re u_k, \quad w_k=(-1)^k\Im u_k.  
\end{equation}

\end{enumerate}

\subsubsection{Asymptotic expansions for large $\tau$}

In order to obtain this expansion, we take the integral representation 
\begin{equation}\label{eq01}
\begin{array}{ll}
\dsp{\Gamma\left(\tfrac12+\mu\right)\sqrt{\pi/2}\,\sinh^\mu\beta \,P_{-\frac12+i\tau}^{-\mu}\left(\cosh\beta\right)
=}
\\
\quad\quad\quad\quad\quad
\dsp{\int_0^\beta\left(\cosh\beta-\cosh t\right)^{\mu-\frac12}\cos \tau t\,dt,\quad \Re\mu>-\tfrac12},
\end{array}
\end{equation}
which is given in \cite[p.~184]{Magnus:1966:FTS}.

The large $\tau$ asymptotics follows from applying the method of stationary phase; see \cite[\S II.3]{Wong:2001:AAI}.

We have the following result
\begin{equation}\label{eq02}
\Gamma\left(\tfrac12+\mu\right)\sqrt{\pi/2}\,\sinh^\mu\beta \,P_{-\frac12+i\tau}^{-\mu}\left(\cosh\beta\right)
\sim
\sum_{n=0}^\infty \frac{A_n}{\tau^{n+1}}+
\sum_{n=0}^\infty \frac{B_n}{\tau^{n+\mu+\frac12}},
\end{equation}
where
\begin{equation}\label{eq03}
A_n=\left.-\sin(n\pi/2)\frac{d^n}{dt^n}\left(\cosh\beta-\cosh t\right)^{\mu-\frac12}\right\vert_{t=0},
\end{equation}
and
\begin{equation}\label{eq04}
B_n=\left.\cos\chi_n\,\frac{\Gamma(n+\mu+\frac12)}{n!}\frac{d^n}{dt^n}\left(\frac{\cosh\beta-\cosh t}{\beta-t}\right)^{\mu-\frac12}\right\vert_{t=\beta},
\end{equation}
where
\begin{equation}\label{eq05}
\chi_n=\tfrac12\left(n-\mu-\tfrac12\right)\pi+\beta \tau.
\end{equation}

The coefficients $A_n$ vanish for even $n$, because of the sine function. They vanish also for odd $n$, because 
in that case the derivatives vanish at $t=0$.

The expansion in (\ref{eq02}), valid for arguments greater than 1 and
large $\tau$, can be written in the form
\begin{equation}\label{eq06}
\begin{array}{lll}
 P_{-\frac12+i\tau}^{m}\left(\cosh\beta\right)&\sim&\sqrt{\Frac{2}{\pi\sinh\beta}}
\cosh(\pi\tau)\Frac{|\Gamma (m+1/2+i\tau)|^2}{\pi}\times\\
&&\sum_{n=0}^\infty \cos\chi_n \left(m+\tfrac12\right)_n\Frac{b_n}{\tau^{n+m+\frac12}},
\end{array}
\end{equation}
where the first few coefficients $b_n$ are given by
\begin{equation}\label{eq07}
\begin{array}{ll}
\dsp{b_0=1,}
\\[8pt]
\dsp{b_1=\frac{(2m-1)x}{4\sinh\beta},}
\\[8pt]
\dsp{b_2=\frac{(2m-1)(-8+(6m-1)x^2)}{96\sinh^2\beta},}
\\[8pt]
\dsp{b_3=\frac{(2m-1)x((-1+4m^2)x^2+16-16m)}{384\sinh^3\beta},}
\end{array}
\end{equation}
and where $\left(m+\tfrac12\right)_n$ is the Pochhammer symbol of $m+\tfrac12$.



 As this expansion holds for small $m$, we will use it for $m=0,1$; then, if larger values
of $m$ are wanted, we will apply forward recursion with the $m$-three term recurrence relation,
as we will later discuss.

\subsection{Computation of $P^m_{-\tfrac12+i\tau}(x)$ for $-1<x<0$} 

Conical functions in the interval $x\in (-1,\,0)$ are computed by means of a stable integral
representation as described in \cite{gil:2009:con}. The starting point is the following integral
representation:

 \begin{equation}\label{Pz1}
P_{-\frac12+i\tau}^{m}(x)=
\frac{\Gamma(\tfrac12+m)\,(1-x^2)^{m/2}\cosh(\pi \tau)}{\pi\sqrt{2\pi}}
\int_{-\infty}^\infty e^{-m\phi(t)}\frac{dt}{\sqrt{x+\cosh t}},
\end{equation}
where
\begin{equation}\label{phitz}
\phi(t)=\ln(x+\cosh t)-i\beta t,\quad \beta= \frac{\tau}{m}.
\end{equation}

In this form, this representation is not suitable for numerical computation because
of the factor $e^{i m \beta t}$ in the integrand: this factor introduces oscillations
which could be very strong for large values of $\tau$. Steepest descent methods \cite[Ch.~5]{Gil:2007:NSF}
can be used for transforming this integral into a stable representation, where oscillations are under control.
In this way, it is possible to obtain the following representation, valid for $-1 <x <0$:

\begin{equation}
\label{Pz2}
\begin{array}{ll}
P_{-\frac12+i\tau}^{m}(x)&=
\Frac{\Gamma(\tfrac12+m)\,(1-x^2)^{m/2}}{\pi\sqrt{2\pi}}
 2\cosh (\pi \tau) e^{-\mu\phi(t_0)} \dsp{\sqrt{\frac{p(1+\beta^2)}{x(p+1)}}}\\
&\displaystyle\int_{0}^\infty e^{-(m+\frac12)\psi_r(s)}\cos((m+\tfrac12)\psi_i(s))\,ds,
\end{array}
\end{equation}
where \begin{equation}\label{betap1}
\beta=\frac{\tau}{m},\quad p=\frac{x}{\sqrt{1+\beta^2(1-x^2)}},
\end{equation}

\begin{equation}\label{phit01}
\phi(t_0)=\ln\frac{x(p+1)}{p(\beta^2+1)}+\beta \arccos \frac{x(1-p\beta^2)}{p(1+\beta^2)}.
\end{equation}
and
\begin{equation}\label{Q5}
\begin{array}{l}
\dsp{\psi_r(s)=\tfrac12\ln\left(1+\frac{4(1+\beta^2)}{1+p}\sigma^2+
\frac{4(1+\beta^2)(1+p^2\beta^2)}{(1+p)^2}\sigma^4\right)},
\\ \\
\dsp{\psi_i(s)=\arctan\frac{\beta(1+p)\sinh s}{1+p+(1-p\beta^2)\sinh^2(s/2)}-\beta s},
\end{array}
\end{equation}
where $\sigma=\sinh(\tfrac12s)$.

\subsection{Three-term recurrence relations}

 Conical functions $P_{-1/2+i\tau}^{m}(x)$ satisfy three-term recurrence relations, which 
are given by:

\begin{equation}
\label{TTRR}
P_{-\frac12+i\tau}^{m+1}(x)+\Frac{2 m x}{\sqrt{1-x^2}}
P_{-\frac12+i\tau}^{m}(x) -
\left((m-\tfrac12)^2+\tau^2\right)P_{-\frac12+i\tau}^{m-1}(x)  =0 
\end{equation}
for $x\in (-1,1)$ and
\begin{equation}
\label{TTRR2}
P_{-\frac12+i\tau}^{m +1}(x)-\Frac{2m x}{\sqrt{x^2-1}}
P_{-\frac12+i\tau}^{m}(x) +
\left((m-\tfrac12)^2+\tau^2\right)P_{-\frac12+i\tau}^{m-1}(x)  =0
\end{equation}
for $x>1$.

These recurrence relations can be used, starting from two initial values, for computing the functions
when they are applied in the direction of stable recursion: either backward or forward.
Also, we will use these relations as a test for the computations.

The stability analysis based on Perron's theorem discussed in \cite{gil:2009:con}, revealed that backward (forward) recursion was 
generally stable for
$x>0$ ($x<0$). However and similarly to other special functions, recurrence relations in the oscillatory regime of the conical functions ($x>x_c>1$)
are not bad conditioned in  both backward and forward directions; then, both recursions are possible. We will use this property
in combination to the asymptotic expansion given in (\ref{eq06})
for computing conical functions for large $\tau$ in the oscillatory regime of the functions.

\section{Overview of the software structure}

The Fortran 90 package includes the main module {\bf Conical}, which includes
the routine {\bf conic}.

In the module {\bf Conical}, the auxiliary 
module {\bf Someconstants} is used. This is a module for the 
computation of the main constants used in 
  the different routines. 
The routines included in {\bf auxil.f90} are also used in the module
{\bf Conical}. Among them, there is a Fortran 90 version of
a Fortran 77 routine for computing the modified Bessel functions $K_{ia}(x)$ and its derivative $K'_{ia}(x)$,
developed by the authors  \cite{Gil:2004:CSM}, \cite{Gil:2004:AMB}.

\section{Description of the individual software components}

The Fortran 90 module {\bf Conical} includes the public routine
{\bf conic} which computes the conical functions $P^m_{-\frac12+i\tau}(x)$, $x>-1$, $m\ge 0$, $\tau >0$.
The calling sequence of this routine is

  \begin{verbatim}
   CALL  conic(x,mu,tau,pmtau,ierr)
  \end{verbatim}

   \noindent
   where the input data are: $x$, $mu$ and $tau$ (arguments of the functions). 
The outputs are the
error flag $ierr$ and the function value $pmtau$.
The possible values of the error flag are: $ierr=0$, successful 
computation; 
$ierr =1$, computation failed due to overflow/underflow; 
$ierr=2$, arguments out of range.  

\section{Testing the algoritm}

\begin{figure}
\caption{Points in the $(x,\tau)$-plane ($x>1$) where the relative error in comparison with the Maple value
in the computation of $P^{95}_{-\frac12+i\tau}(x)$
is $\sim 10^{-10}$. At the rest of tested points in the $(x,\tau)$-plane,
the accuracy was found $\sim 10^{-12}$ or better.  The curve $y=\sqrt{1+\beta^2}/\beta$, where $\beta= \tau/95$,
and the regions where different approximations are used in the algorithm for $x>1$, are also shown in the figure.
\label{Fig3}}
\begin{center}
\epsfxsize=13cm \epsfbox{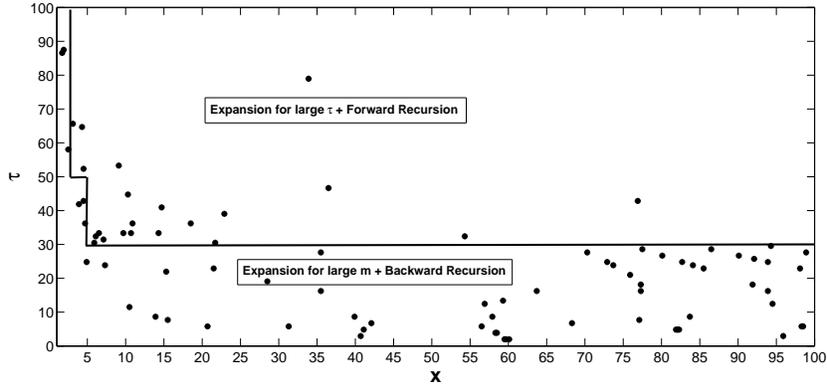}
\end{center}
\end{figure}

 Two kind of tests have been considered for testing the accuracy of the computed values of the 
conical functions: direct comparison against Maple and a single step
of the three-term recurrence relations
(\ref{TTRR}) and (\ref{TTRR2}). At most of the points of the tested parameter space $(x,\tau,m)$ the 
comparison against Maple shows that the accuracy was $\sim 10^{-12}$ or better ($10^{-13}-10^{-14}$
in $-1<x<0$, included in the monotonic region). It is important
to point out that for $x>1$, asymptotic expansions using the modified Bessel functions $K_{ia}(x)$
and $K'_{ia}(x)$ (\ref{PzwK}) are considered and that the accuracy in the computation of these 
functions is $\sim 5\times 10^{-13}$, as explained in \cite{Gil:2004:AMB}. So, the accuracy
for computing these functions limit the attainable accuracy in the computation
of conical functions. 
On the other hand, at some points
of the oscillatory region of the conical functions, the tested accuracy was $\sim 10^{-10}$. 
This is apparent in Figure (\ref{Fig3}), where points in the $(x,\tau)$-plane with a relative error (in comparison with the Maple value)
$\sim 10^{-10}$ in the computation of $P^{95}_{-\frac12+i\tau}(x)$, are plotted. The figure also shows 
the curve $y=\sqrt{1+\beta^2}/\beta$, where $\beta=\tau/95$. This curve
is the frontier between the monotonic and the oscillatory regions for the conical function $P^{95}_{-\frac12+i\tau}(x)$.
Additionally, the approximations used in the algorithm for $x>1$ are indicated in the figure. As can be seen, the density
of plotted points is larger in the region where asymptotic expansions in terms of modified Bessel functions are used. 
Finally, it is important also to note that in the oscillatory region the zeros of the conical functions are found and at these
points relative error losses its meaning.

\section{Test run description}

The Fortran 90 test program {\bf testcon.f90} includes the computation
of 25 function values and their comparison with the corresponding
pre-computed results. Also, a single step of the three-term recurrence relations
(\ref{TTRR}) and (\ref{TTRR2}) is tested for several values of the parameters $(x,\,\tau,\,m)$.

\section{Acknowledgements}

The authors thank the referee for useful comments.
The authors acknowledge financial support from 
{\emph{Ministerio de Ciencia e Innovaci\'on}}, project MTM2009-11686. NMT acknowledges financial support 
from {\emph{Gobierno of Navarra}, Res. 07/05/2008}.








\end{document}